\documentclass[10pt, a4paper]{amsart}
\usepackage[dvips]{epsfig}
\usepackage{amsmath}
\usepackage{amssymb}
\usepackage{amsfonts}
\usepackage{amsthm}
\usepackage{amsbsy}
\usepackage{amsgen}
\usepackage{amscd}
\usepackage{amsopn}
\usepackage{amstext}
\usepackage{amsxtra}
\usepackage{xypic}

\newtheorem{theorem}{Theorem}[section]

\newtheorem{defn}[theorem]{Definition}

\theoremstyle{definition}

\numberwithin{equation}{section}

\begin{document}

\title[A homotopy theory for enrichment in simplicial modules]
{A homotopy theory for enrichment in simplicial modules}
\author{Alexandru E. Stanculescu}
\address{Department of Mathematics and Statistics, McGill University, 805 Sherbrooke Str. West, Montr\'eal,
Qu\'ebec, Canada, H3A 2K6} \email{stanculescu@math.mcgill.ca}

\date{December 6, 2007}

\begin{abstract}
We put a Quillen model structure on the category of small categories
enriched in simplicial $k$-modules and non-negatively graded chain complexes
of $k$-modules, where $k$ is a commutative ring.
The model structure is obtained by transfer from the model structure
on simplicial categories due to J. Bergner.
\end{abstract}
\maketitle
\section{Introduction: DK-equivalences and DK-fibrations}
1.1. Let $\mathbf{Cat}$ the category of small categories. It has a
\emph{natural} model structure in which a cofibration is a functor
monic on objects, a weak equivalence is an equivalence of categories
and a fibration is an isofibration \cite{JT0}. The fibration weak
equivalences are the equivalences surjective on objects.

Let $\mathcal{V}$ be a monoidal model category \cite{SS2} with unit
$I$. We denote by $\mathfrak{W}$ the class of weak equivalences of
$\mathcal{V}$, by $\mathfrak{Fib}$ the class of fibrations and by
$\mathfrak{Cof}$ the class of cofibrations.

The small $\mathcal{V}$-categories together with the
$\mathcal{V}$-functors between them form a category written
$\mathcal{V}\mathbf{Cat}$. Let $\mathcal{M}$ be a class of maps of
$\mathcal{V}$. We say that a $\mathcal{V}$-functor
$f:\mathcal{A}\rightarrow \mathcal{B}$ is \emph{locally in}
$\mathcal{M}$ if for each pair $x,y\in \mathcal{A}$ of objects, the
map $f_{x,y}:\mathcal{A}\rightarrow \mathcal{B}$ is in
$\mathcal{M}$.

We have a functor $[\_]_{\mathcal{V}}:\mathcal{V}\mathbf{Cat}
\rightarrow \mathbf{Cat}$ obtained by change of base along the
(symmetric monoidal) composite functor

\[
   \xymatrix{
   \mathcal{V} \ar[r]^{\gamma} & Ho(\mathcal{V})
   \ar[rr]^{Hom_{Ho(\mathcal{V})}(I,\_)} & & Set.\\
}
\]

\begin{defn} Let $f:\mathcal{A}\rightarrow \mathcal{B}$ be a
morphism in $\mathcal{V}\mathbf{Cat}$.

1. The morphism $f$ is $\mathrm{homotopy \ essentially \
surjective}$ if the induced functor
$[f]_{\mathcal{V}}:[\mathcal{A}]_{\mathcal{V}}\rightarrow
[\mathcal{B}]_{\mathcal{V}}$ is essentially surjective.

2. The morphism $f$ is a $\mathrm{DK-\ equivalence}$ if it is
homotopy essentially surjective and locally in $\mathfrak{W}$.

3. The morphism $f$ is a $\mathrm{DK-\ fibration}$ if it satisfies
the following two conditions.

(a) $f$ is locally in $\mathfrak{Fib}$.

(b) For any $x\in \mathcal{A}$, and any isomorphism
$v:[f]_{\mathcal{V}}(x)\rightarrow y'$ in
$[\mathcal{B}]_{\mathcal{V}}$, there exists an isomorphism
$u:x\rightarrow y$ in $[\mathcal{A}]_{\mathcal{V}}$ such that
$[f]_{\mathcal{V}}(u)=v$. That is, if $[f]_{\mathcal{V}}$ is an
isofibration.
\end{defn}
One can easily see that a morphism $f$ is a DK-equivalence and a
DK-fibration iff $f$ is surjective on objects and locally in
$\mathfrak{W}\cap \mathfrak{Fib}$. The class of maps having the left
lifting property with respect to the $\mathcal{V}$-functors
surjective on object and locally in $\mathfrak{W}\cap
\mathfrak{Fib}$ is generated by the map $u:\emptyset \rightarrow
\mathcal{I}$, where $\mathcal{I}$ is the $\mathcal{V}$-category with
a single object $\ast$ and $\mathcal{I}(\ast,\ast)=I$, together with
the maps $$\bar{2}_{i}:\bar{2}_{A}\rightarrow \bar{2}_{B},$$ where
$i$ is a generating cofibration of $\mathcal{V}$. Here the
$\mathcal{V}$-category $\bar{2}_{A}$ has objects 0  and 1, with
$\bar{2}_{A}(0,0)=\bar{2}_{A}(1,1)=I$, $\bar{2}_{A}(0,1)=A$ and
$\bar{2}_{A}(1,0)=\emptyset$.

 1.2. Let $k$ be a commutative ring. We denote by
$\mathbf{SMod}_{k}$ the category of simplicial $k$-modules and by
$Ch^{+}(k)$ the category of non-negatively graded chain complexes of
$k$-modules. The purpose of this note is to prove the following
theorem.
\begin{theorem}
Let $\mathcal{V}$ be one of the categories $\mathbf{SMod}_{k}$ or
$Ch^{+}(k)$. Then $\mathcal{V}\mathbf{Cat}$ admits a model structure
in which the weak equivalences are the DK-equivalences and the
fibrations are the DK-fibrations.
\end{theorem}
To prove this result we use the (similar) model structure on
simplicial categories \cite{Be} and Quillen's path object argument
(\cite{SS1}, Lemma 2.3(2) and \cite{BM}, 2.6). An explicit
description of a cofibration of $\mathcal{V}\mathbf{Cat}$ can be
given \cite{St}.

1.3. The proof of theorem 1.2 relies decisively on the construction
of path objects for dg-categories due to G. Tabuada (\cite{Ta2},
4.1). In fact, our attempt to understand his construction led us to
the proof of our result.

1.4. In \cite{To}, B. To\"{e}n characterised the maps in the
homotopy category of dg-categories, where the category of
dg-categories has a model structure in which the weak equivalences
are the DK-equivalences and the fibrations are the DK-fibrations.
One can show that his results (\emph{loc. cit.}, Thm. 4.2 and 6.1)
hold for $\mathcal{V}\mathbf{Cat}$, where $\mathcal{V}$ is
$\mathbf{SMod}_{k}$ or $Ch^{+}(k)$.

\emph{Note added in proof.} After the completion of this work we
learned about the existence of a paper by G. Tabuada \cite{Ta1},
which treats the same subject matter, and more, but differently. One
can see that the model structure proposed in theorem 1.2 coincides
with the one in \cite{Ta1}, although the classes of fibrations and
cofibrations are not explicitly identified in \emph{loc. cit.} On
the other hand, Tabuada shows that the model structures on
$\mathbf{SMod}_{k}\mathbf{Cat}$ and $Ch^{+}(k)\mathbf{Cat}$ are
Quillen equivalent, an issue that we have initially neglected. One
can easily give a proof of this fact, adapted to our context, using
section 2.2 below and the general results of \cite{St}.

\section{Categories enriched in $\mathbf{SMod}_{k}$ and $Ch^{+}(k)$}
2.1. The category $\mathbf{SMod}_{k}$ is a closed symmetric monoidal
category with tensor product defined pointwise and unit $ck$, where
$(ck)_{n}=k$ for all $n\geq 0$. A model structure on
$\mathbf{SMod}_{k}$ is obtained by transfer from the category
$\mathbf{S}$ of simplicial sets, regarded as having the classical
model structure, via the free-forgetful adjunction
$$k:\mathbf{S} \rightleftarrows \mathbf{SMod}_{k}:U.$$
All objects are fibrant and the model structure is simplicial. The
functor $k$ is strong symmetric monoidal (and it preserves the
unit), hence $\mathbf{SMod}_{k}$ is a monoidal model category. The
adjunction ($k,U$) induces an adjunction
$$k':\mathbf{SCat} \rightleftarrows \mathbf{SMod}_{k}\mathbf{Cat}:U'.$$
We claim that a map $f$ of $\mathbf{SMod}_{k}\mathbf{Cat}$ is a
DK-equivalence (resp. DK-fibration) iff $U'(f)$ is a weak
equivalence (resp. fibration) in the Bergner model structure on
$\mathbf{SCat}$ \cite{Be}. Clearly, $f$ is locally in $\mathfrak{W}$
(resp. $\mathfrak{Fib}$) iff $U'(f)$ is locally in $\mathfrak{W}$
(resp. $\mathfrak{Fib}$). In the induced adjoint pair
$$Lk:Ho(\mathbf{S}) \rightleftarrows Ho(\mathbf{SMod}_{k}):RU,$$
the functor $Lk$ is strong symmetric monoidal and  preserves the unit object, hence one has
a natural isomorphism of functors
$$\eta:[\_]_{\mathbf{SMod}_{k}}\cong [\_]_{\mathbf{S}}U':\mathbf{SMod}_{k}\mathbf{Cat}\rightarrow \mathbf{Cat}$$
such that for all $\mathcal{A} \in \mathbf{SMod}_{k}\mathbf{Cat}$, $\eta_{\mathcal{A}}$ is the identity on objects.
The rest of the claim follows from this observation.

2.2. Consider the normalized chain complex functor $N:\mathbf{SMod}_{k} \rightarrow Ch^{+}(k)$. It was shown in
(\cite{SS2}, 4.3) that $N$ is part of a weak monoidal Quillen equivalence
$$N:\mathbf{SMod}_{k}\rightleftarrows Ch^{+}(k):\Gamma$$ in which both functors preserve the unit objects.
Therefore the composite adjunction
$$Nk:\mathbf{S}\rightleftarrows Ch^{+}(k):U\Gamma$$ is a weak monoidal Quillen
pair with $Nk$ preserving the unit object. The functor $U\Gamma$
induces a functor $(U\Gamma)':Ch^{+}(k)\mathbf{Cat}\rightarrow
\mathbf{SCat}$ which has a left adjoint $F$ defined "fibrewise". We
claim that a map $f$ in $Ch^{+}(k)\mathbf{Cat}$ is a DK-equivalence
(resp. DK-fibration) iff $(U\Gamma)'(f)$ is a weak equivalence
(resp. fibration) in the Bergner model structure on $\mathbf{SCat}$.
For this, it is enough to remark that in the induced adjunction
$$L(Nk):Ho(\mathbf{S}) \rightleftarrows Ho(Ch^{+}(k)):R(U\Gamma),$$
the functor $L(Nk)$ is strong monoidal and preserves the unit object,
and then conclude as in 2.1.

2.3. In order to prove theorem 1.2, it suffices to apply Quillen's path object argument
to the adjunctions $(k',U')$ and $(F,(U\Gamma)')$. This will be achieved in the next section.

\section{Cocategory object structure on the interval}
3.1. Let $\mathcal{V}$ be a monoidal model category with cofibrant
unit $I$ and all objects fibrant. We write $Y^{X}$ for the internal
hom of two objects $X,Y$ of $\mathcal{V}$. We say that $\mathcal{V}$
has a \emph{cocategory interval} if there is a cocategory object
structure
\[
   \xymatrix{
  I \ar @<-4pt> [r]_{d_{1}} \ar @<10pt> [r]^{d_{0}} & I[1] \ar @<-1pt> [l]_{p}
  \ar@<10pt> [r]^{i_{0}} \ar@<2pt> [r]^{c} \ar@<-4pt> [r]_{i_{1}} & I[2]\\
  }
  \]
such that
\[
   \xymatrix{
  I\sqcup I \ar[rr]^{\bigtriangledown} \ar[dr]_{d_{0}\sqcup d_{1}} & & I\\
  & I[1] \ar[ur]_{p}\\
  }
  \]
is a cylinder object for $I$. The map $c$ denotes the cocomposition.

\emph{Examples}. $(a)$ The standard example is when
$\mathcal{V}=\mathbf{Cat}$ as in 1.1. Here $I[1]$ is the
"free-living" isomorphism and $I[2]$ is the groupoid with three
objects and one isomorphism between any two objects. We leave to the
reader the task to identify all the maps involved.

$(b)$ The case which interests us is when $\mathcal{V}=Ch^{+}(k)$. The interval $I[1]$ is well known to be
$...\rightarrow 0\rightarrow ke\overset{\partial}\rightarrow ka\oplus kb$, where $\partial(e)=b-a$.
 The maps $d_{0}$ and $d_{1}$ are the inclusions, and the map $p$ is $a,b \mapsto 1$. The object $I[2]$ is
$$...\rightarrow 0\rightarrow ke_{1}\oplus ke_{2}\overset{\partial} \rightarrow ka_{0}\oplus ka_{1} \oplus ka_{2},$$
where $\partial(e_{1})=a_{1}-a_{0}$ and $\partial(e_{2})=a_{2}-a_{1}$. The cocomposition $c$ is given
by $e\mapsto e_{1}+e_{2}$, $a\mapsto a_{0}$ and $b\mapsto a_{2}$. The map $i_{0}$ is given by
$e\mapsto e_{1}$, $a\mapsto a_{0}$ and $b\mapsto a_{1}$; the map $i_{1}$ is given by
$e\mapsto e_{2}$, $a\mapsto a_{1}$ and $b\mapsto a_{2}$.

$(c)$ Since the functor $\Gamma$ from 2.2 preserves the unit object
and is an equivalence of categories, we obtain that
$\mathbf{SMod}_{k}$ has a cocategory interval.

$(d)$ Let $k$ be a field and let $H$ be a finite dimensional cocommutative Hopf algebra over $k$.
We let $\mathcal{V}=_{H}Mod$, the category of left $H$-modules, and we view $\mathcal{V}$ as having the
stable model structure \cite{Ho}. Let $u$ (resp. $\epsilon$) be the unit (resp. counit) of $H$. A cylinder object
for $k$ is
\[
   \xymatrix{
  k\oplus k \ar[rr]^{\bigtriangledown} \ar[dr]_{id\oplus u} & & k\\
  & k\oplus H \ar[ur]_{p:=(id,\epsilon)}\\
  }
  \]
The maps $d_{0}$ and $d_{1}$ are given by $d_{0}(1)=(1,0)$ and
$d_{1}(1)=(0,1_{H})$. We set $I[2]=k\oplus k\oplus H$ and
$c(\alpha,h)=(\alpha,0,h)$. The map $i_{0}$ is $(\alpha,h)\mapsto
(\alpha,\epsilon(h),0)$ and the map $i_{1}$ is $(\alpha,h)\mapsto
(0,\alpha,h)$. One can check that the resulting gadget is a
cocategory interval with $I[1]=k\oplus H$. It is easy to see that
$I[1]$ is an "interval with a coassociative and cocommutative
comultiplication" in the sense of (\cite{BM}, page 813).

3.2. We shall now construct (DK-)path objects for
$\mathcal{V}$-categories, where $\mathcal{V}$ is as in 3.1. In the
case of dg-categories, the construction is due to G. Tabuada
(\cite{Ta2}, 4.1). Let $\mathcal{A}\in \mathcal{V}\mathbf{Cat}$. We
first construct a factorisation
\[
   \xymatrix{
  \mathcal{A} \ar[rr]^{\bigtriangleup} \ar[dr]_{i_{0}} & & \mathcal{A} \times \mathcal{A}\\
  & P_{0}\mathcal{A} \ar[ur]_{(s,t)}\\
  }
  \]
such that $i_{0}$ is locally in $\mathfrak{W}$ and $(s,t)$ is locally in $\mathfrak{Fib}$. An object of
$P_{0}\mathcal{A}$ is a map $f:a\rightarrow b$ of $[\mathcal{A}]_{\mathcal{V}}$. If $f_{0}:a_{0}\rightarrow b_{0}$ and
$f_{1}:a_{1}\rightarrow b_{1}$ are two objects of $P_{0}\mathcal{A}$, we define $P_{0}\mathcal{A}(f_{0},f_{1})$
to be the limit of the diagram
\[
   \xymatrix{
  \mathcal{A}(a_{0},a_{1}) \ar[dr]_{f_{1*}} & & \mathcal{A}(a_{0},b_{1})^{I[1]} \ar[dl]^{t} \ar[dr]_{s}
   & & \mathcal{A}(b_{0},b_{1}) \ar[dl]^{f_{0}^{*}}\\
  & \mathcal{A}(a_{0},b_{1}) & & \mathcal{A}(a_{0},b_{1})\\
  }
  \]
The unit of $P_{0}\mathcal{A}(f,f)$ is induced by the adjoint transpose of
$I[1]\overset{p}\rightarrow I \overset{f}\rightarrow \mathcal{A}(a,b)$. Let $f_{i}:a_{i}\rightarrow b_{i}$
$(i=0,2)$ be three objects of $P_{0}\mathcal{A}$ and let $A_{i}=P_{0}\mathcal{A}(f_{i},f_{i+1})$ $(i=0,1)$.
We denote by $p_{i}$ (resp. $q_{i}$) the canonical map $A_{i}\rightarrow \mathcal{A}(a_{i},a_{i+1})$
(resp. $A_{i}\rightarrow \mathcal{A}(b_{i},b_{i+1})$) $(i=0,1)$. The pair $(p_{i},q_{i})$ gives rise to a commutative diagram
\[
   \xymatrix{
 A_{i} \ar[r]^{p_{i}} \ar[d]_{j_{1,i}} & \mathcal{A}(a_{i},a_{i+1}) \ar[d]^{f_{i+1 \ast}}\\
 I[1]\otimes A_{i} \ar[r]^{H_{i}} & \mathcal{A}(a_{i},b_{i+1})\\
 A_{i} \ar[u]^{j_{0,i}} \ar[r]^{q_{i}} & \mathcal{A}(b_{i},b_{i+1}) \ar[u]_{f_{i}^{\ast}}\\
 }
  \]
where $j_{k,i}=d_{k}\otimes A_{i}$ $(k=0,1)$. Observe that in order to define a map
$A_{0}\otimes A_{1} \rightarrow P_{0}\mathcal{A}(f_{0},f_{2})$
it suffices to find a map $G:A_{0}\otimes A_{1} \rightarrow \mathcal{A}(a_{0},b_{2})^{I[1]}$
which makes commutative the diagram
\[
   \xymatrix{
   & \mathcal{A}(a_{0},b_{2})\\
 A_{0}\otimes A_{1} \ar[r]^{G} \ar[dr]_{f_{2\ast}(p_{0}\otimes p_{1})} \ar[ur]^{f_{0}^{\ast}(q_{0}\otimes q_{1})}
  & \mathcal{A}(a_{0},b_{2})^{I[1]} \ar[d]^{t} \ar[u]_{s}\\
 & \mathcal{A}(a_{0},b_{2})\\
 }
  \]
We define the map $G_{1}$ as the composite
$$I[1]\otimes A_{0}\otimes A_{1} \overset{H_{0}\otimes A_{1}}\longrightarrow \mathcal{A}(a_{0},b_{1})\otimes A_{1}
\overset{id\otimes q_{1}}\longrightarrow \mathcal{A}(a_{0},b_{1}) \otimes \mathcal{A}(b_{1},b_{2})\rightarrow \mathcal{A}(a_{0},b_{2}).$$
Then $G_{1}$ is a "homotopy" between $f_{0}^{\ast}(q_{0}\otimes q_{1})$ and $p_{0}\otimes q_{1}$.
We define the map $G_{2}$ as the composite
$$I[1]\otimes A_{0}\otimes A_{1} \overset{A_{0}\otimes H_{1}}\longrightarrow A_{0} \otimes \mathcal{A}(a_{1},b_{2})
\overset{p_{0}\otimes id}\longrightarrow \mathcal{A}(a_{0},a_{1}) \otimes \mathcal{A}(a_{1},b_{2})\rightarrow \mathcal{A}(a_{0},b_{2}).$$
Then $G_{2}$ is a "homotopy" between $p_{0}\otimes q_{1}$ and $f_{2\ast}(p_{0}\otimes p_{1})$. The two homotopies
induce a map $$A_{0}\otimes A_{1}\rightarrow \mathcal{A}(a_{0},b_{2})^{I[1]} \times_{\mathcal{A}(a_{0},b_{2})}
\mathcal{A}(a_{0},b_{2})^{I[1]}$$ such that the diagram
\[
   \xymatrix{
   & \mathcal{A}(a_{0},b_{2})\\
 A_{0}\otimes A_{1} \ar[r] \ar[dr]_{\bar{G_{2}}} \ar[ur]^{\bar{G_{1}}}
  & \mathcal{A}(a_{0},b_{2})^{I[1]}\times_{\mathcal{A}(a_{0},b_{2})}
\mathcal{A}(a_{0},b_{2})^{I[1]} \ar[d] \ar[u]\\
 & \mathcal{A}(a_{0},b_{2})\\
 }
  \]
commutes, where $\bar{G_{i}}$ is the adjoint transpose of $G_{i}$ $(i=0,1)$. Since
$\mathcal{A}(a_{0},b_{2})^{I[1]}$ is a category object in $\mathcal{V}$, we have then a map
$$A_{0}\otimes A_{1}\rightarrow \mathcal{A}(a_{0},b_{2})^{I[1]} \times_{\mathcal{A}(a_{0},b_{2})}
\mathcal{A}(a_{0},b_{2})^{I[1]} \overset{m} \longrightarrow \mathcal{A}(a_{0},b_{2})^{I[1]}$$
which is the required map $G$.

In this way $P_{0}\mathcal{A}$ becomes a $\mathcal{V}$-category. The association
 $i_{0}:Ob(\mathcal{A})\rightarrow Ob(P_{0}\mathcal{A})$, $a\mapsto (id_{a}:a\rightarrow a)$,
$(i_{0})_{a,b}=\mathcal{A}(a,b)^{p}$, is a $\mathcal{V}$-functor
$\mathcal{A}\rightarrow P_{0}\mathcal{A}$. By construction, the maps
$s,t:P_{0}\mathcal{A}\rightarrow \mathcal{A}$,
$s(f_{0}:a_{0}\rightarrow b_{0})=a_{0}$, $t(f_{0}:a_{0}\rightarrow
b_{0})=b_{0}$, $s_{f_{0},f_{1}}=p_{0}$ and $t_{f_{0},f_{1}}=q_{0}$,
are $\mathcal{V}$-functors. One clearly has
$(s,t)i_{0}=\bigtriangleup$. Moreover, $(s,t)$ is locally
 in $\mathfrak{Fib}$ since $(s,t)_{f_{0},f_{1}}:A_{0}\rightarrow \mathcal{A}(a_{0},a_{1}) \times \mathcal{A}(b_{0},b_{1})$ is the pullback
\[
   \xymatrix{
 & \mathcal{A}(a_{0},a_{1})\times \mathcal{A}(b_{0},b_{1}) \ar[d]^{f_{1\ast}\times f_{0}^{\ast}}\\
\mathcal{A}(a_{0},b_{1})^{I[1]} \ar[r]^{(s,t)} & \mathcal{A}(a_{0},b_{1})\times \mathcal{A}(a_{0},b_{1}).\\
}
  \]
Next, let $P\mathcal{A}$ be the full sub-$\mathcal{V}$-category of $P_{0}\mathcal{A}$ whose
objects consist of isomorphisms $f:a\rightarrow b$ of $[\mathcal{A}]_{\mathcal{V}}$. Then $i_{0}$
factors through $P\mathcal{A}$. The resulting factorisation
\[
   \xymatrix{
  \mathcal{A} \ar[rr]^{\bigtriangleup} \ar[dr]_{i} & & \mathcal{A} \times \mathcal{A}\\
  & P\mathcal{A} \ar[ur]_{(s,t)}\\
  }
  \]
is the desired (DK-)path object: a lengthy but straightforward
computation shows that $i$ is homotopy essentially surjective and
that $[(s,t)]_{\mathcal{V}}$ is an isofibration.
\\

\textbf{Acknowledgments.} I would like to thank Professor Andr\'{e}
Joyal for many useful discussions about the subject and for
pronouncing the word "cocategory".

\end{document}